\documentclass[12pt,a4paper]{article}
\usepackage{amsmath,amsfonts,amssymb,amsthm,}
\usepackage{color}
\usepackage[english]{babel}
\usepackage{vmargin}

\setmarginsrb{23mm}{20mm}{23mm}{20mm}{10mm}{10mm}{10mm}{10mm}

\newtheorem{theo}{Theorem}
\newtheorem{coro}{Corollary}
\newtheorem{propo}{Proposition}

\theoremstyle{definition}

\newtheorem*{Remark}{Remark}

\theoremstyle{remark}

\newcommand{\C}{\mathbb{C}}
\newcommand{\D}{\mathbb{D}}

\newcommand{\T}{\mathbb{T}}

\newcommand{\diam}{\operatorname{diam}}

\definecolor{blau}{rgb}{0.1,0.0,0.9}
\definecolor{violet}{rgb}{0.54, 0.17, 0.89}
\newcommand{\blue}{\color{blau}}

\newcommand{\kom}[1]{}
\renewcommand{\kom}[1]{{\bf \blue /#1/}}

\newcounter{komcounter}
\numberwithin{komcounter}{section}

\title{Quasiconformal properties of $Q_{p,0}$ curves and Dirichlet-type curves}

\author{Mar\'ia\ J.\ Gonz\'alez}

\date{}

\begin{document}
	
	\maketitle
	
	\baselineskip=6mm
	\parskip=2.5mm
	\let\thefootnote\relax
	\let\thefootnote\relax

	\begin{abstract}
	Let $\Gamma$ be a closed Jordan curve, and  $f$ the conformal mapping that sends the unit disc $\D$  onto the interior  domain of $\Gamma$. If $\log f'$ belongs to the  Dirichlet space $\mathcal{D}$, we call $\Gamma$ a Weil-Petersson curve. The purpose of this note is to  extend recent results, obtained by   G. Cui and  Ch. Bishop in the case of Weil-Petersson curves, to the case when  $\log f'$ belongs to either some $Q_{p,0},$ space, for  $0<p\leq 1$, or to some weighted-Dirichlet space contained in $\mathcal{D}$.  More precisely, we will characterize the quasiconformal extensions of $f$, and describe some  of the geometric properties of $\Gamma$, that arise in this context.

	\end{abstract}

	\footnotetext{The author was supported by grant no. MTM2017-85666-P, Spain}.

	\section*{Introduction}


	 A closed Jordan  curve $\Gamma$ is called a Weil-Petersson curve if the conformal map $f$ from the unit disc $\D$ onto the interior domain of $\Gamma$ satisfies that $\log{f'}$ belongs to the Dirichlet space $\mathcal{D}$, that is, 
	 $$
	 \int_{\D}| (\log f')'|^2 dA(z)=\int_{\D} \bigg\vert\frac{f''(z)}{f'(z)~}\bigg\vert^2 dA(z)<\infty.
	 $$
	 
	 Weil-Petersson curves have drawn much attention recently, we refer the reader to the monograph by Ch. Bishop \cite{b} and the references therein, where numerous  characterizations of these type of curves are provided in different settings, including geometric function theory and Teichm\"uller theory.\\
	 Motivated by some of the questions mentioned in this monograph, we consider in this note conformal mappings $f$ on the unit disc such that $\log f'$ belongs to certain spaces of analytic functions  closely related to the Dirichlet space, namely: Dirichlet-type spaces that we will denote by $\mathcal{D}_{\log_p}, p\geq 0$, where we integrate against a weight that blows up logarithmically when we approach the boundary of the unit disc, and the  $Q_{p,0} $ spaces, $0<p\leq 1$.

	 We say that an analytic function $g$ defined in $\D$ is in $ \mathcal{D}_{\log_p}, ~p\geq 0$ if $$\int_\D |g'|^2 \left(\log \frac{1}{1-|z|}\right)^p<\infty.$$ These spaces are contained in the Dirichlet space $\mathcal{D}$, which corresponds to $p=0$.
	 
	 An analytic function $g$ on $\D$ is said to be in  $Q_{p,0},~p>0$ if the measure $d\mu=|g'(z)|^2 (1-|z|^2)^p dA(z)$ is a $p$-vanishing Carleson measure  on $\D$, that is,
	 $$
	 \lim_{|I|\to 0}\frac{1}{|I|^p}\int_{Q_I} |g'(z)|^2 (1-|z|^2)^p dA(z)= 0,
	 $$
	 where $I$ is an arc on $\T=\partial \D$, and $Q_I$ is the corresponding Carleson box $Q_I=\{z\in \D;~ 1-|I|\leq |z|<1, z/|z|\in I\}\subset \D $.\\
	 \indent
	 If $p > 1$, the space  $Q_{p,0}$ coincides with the little Bloch space $\mathcal{B}_0$, and for $p=1, Q_{1,0} = VMOA$, the space of analytic functions on $ \D$
	 with vanishing mean oscillation. In fact, for $0<p<1, ~ \mathcal{D} \subset Q_{p,0} \subset VMOA\subset \mathcal{B}_0 $.
	 
	 Our first results are related to  quasiconformal extensions of the conformal mappings. Recall that a global homeomorphism  on the plane $~\rho$ is called quasiconformal if it   preserves orientation, belongs to the Sobolev class $W^{1,2}_{loc}(\mathbb{C})$, and  satisfies the Beltrami equation
	 $\overline{\partial}\rho-\mu \partial\rho=0$, where $\mu$ is a measurable function, called the complex dilatation, such that $\|\mu\|_{\infty}<1$.
	 
	 \begin{theo}\label{theo 1} Let $f:\D \to \Omega$ be a conformal map of the unit disc $\D$ onto a bounded domain $\Omega$ whose boundary is a Jordan curve $\Gamma$. For any $p\geq 0$, if	$\log{f'}\in \mathcal{D}_{\log_p},$ then $f$ can be extended to a global quasiconformal mapping  with complex dilatation $\mu$ such that  
	 	\begin{equation}\label{D}
	 	\int_{\C\backslash \D }~ \frac{|\mu(z)|^2}{(|z|^2-1)^2 }\left(\log \frac{|z|}{|z|-1}\right)^p~ dA(z)<\infty.
	 	\end{equation}
	 	On the other hand, if there exists a quasiconformal extension of $f$ to the whole plane with complex dilatation $\mu$ satisfying (\ref{D}), and such that $||\mu||_\infty <1/2$ in a neighborhood of the boundary of the unit disc, then $\log{f'}\in \mathcal{D}_{\log_p}.$ 
	 	
	 \end{theo}
 \begin{Remark}  Note that if $\log{f'}\in \mathcal{D}_{\log_p},$ for $p>2$, then $f'$ is Dini smooth in $\overline\D$. In particular, $\Gamma$ is a Dini-smooth curve. For the sake of completeness, we will give the proof of this remark in the next section.
 \end{Remark}

 In the context of $Q_{p,0}$ spaces, we obtain the following: 
 
 \begin{theo} \label {Theorem 2}Let $f:\D \to \Omega$ be the conformal map of the unit disc $\D$ onto a bounded domain $\Omega$  whose boundary is a Jordan curve $\Gamma$. If	 $\log{f'}\in Q_{p,0},~0<p\leq 1$, then  $d\nu=\frac{|\mu(z)|^2}{(|z|-1)^{2-p} }~ dA(z)$ is a p-vanishing Carleson measure in $\C\backslash \D$, that is,  for any arc $I\in \T$, 
 	\begin{equation}
 	\lim_{I\to 0} \frac{1}{|I|^p}	\int_{\tilde{Q}_I } \frac{|\mu(z)|^2}{(|z|-1)^{2-p} }~ dA(z)=0,\label{con}
 	\end{equation}
 	where $\tilde{Q}_I$ is the Carleson box ~$\tilde{Q}_I=\{z\in \C\backslash \D;~ 1< |z|\leq 1+|I|, z/|z|\in I\}\subset\C\backslash \D$ associated to $I$.\\
 	Conversely, if there exists a quasiconformal extension of $f$ to the whole plane with complex dilatation $\mu$ satisfying (\ref{con}), and such that $||\mu||_\infty<(1+p)/2$ in a neighborhood of the boundary of the unit disc, then  $\log{f'}\in Q_{p,0}.$
 \end{theo}

Let us just point out that the  assumptions on the $||\mu||_\infty$ norm are rather technical, nevertheless they are needed in our proofs. It is is very likely that the conclusions  still hold without this extra assumption. 
	 
	Particular instances of these results are already known, namely  the Dirichlet space, that corresponds to $p=0$ in Theorem 1, and  the $VMOA$ space that corresponds to $p=1$ in Theorem 2. The Dirichlet result is  due to G. Cui  (\cite{c}), the proof involves the  Schwarzian derivative of $f$, and in fact, the result is obtained without using  the extra assuption $||\mu||_\infty < 1/2$. Although we will follow a different   approach, let us  mention here that  characterizations of the Schwarzian derivative of $f$ when $\log f'\in Q_p$ or $Q_{p,0}$ can be found in \cite{PP} and \cite{PR}. Our  proofs will follow closely the argument presented by E. Dyn'kin in (\cite{d}, sect 4) where the $VMOA$ case is established.\\
	\indent
In order to contextualize these results, let us consider the following extension of the conformal mapping $f$, that we also denote by $f$:
\begin{equation}\label{beta}
f(z) = f(1/\bar{z}) + f'(1/\bar{z})(z-1/\bar{z}), \quad  |z| > 1.
\end{equation}
\noindent
An immediate computation shows that,
$$
\mu(z)=\frac{\partial_{\bar{z}} f(z)}{\partial_z f(z)}= -\frac{1}{\bar{z}^2}\left(z- \frac{1}{\bar{z}}\right) \frac{f''(1/\bar{z})}{f'(1/\bar{z})},  \quad  |z| > 1.
$$
\noindent
Hence, if we define 
$$
\eta(z) =(1 - |z|)~\bigg\vert\frac{f''(z)}{f'(z)}\bigg\vert, \quad  z\in \D,
$$
it holds that for  $|z|> 1$, 
\begin{equation}\label{2}
|\mu(z)|= \frac{|z|+1}{|z|^2}~ \eta(1/\bar{z}).
\end{equation}

In general, the mapping (\ref{beta}) is not homeomorphic, but Becker and Pommerenke \cite{bp} proved that (\ref{beta}) is indeed a quasiconformal extension
of $f$ to a neighborhood of $\T=\partial \D$ if $f(\D)$ is a Jordan domain and $\limsup_{|z|{\rightarrow 1-}}\eta(z) < 1$. Thus, in view of (\ref{2}), the statements in Theorem 1 and Theorem 2 are not surprising. \\
\indent 
The corresponding result when	$\log{f'}\in \mathcal{B}_0$, i.e. $ \eta(z)= (1-|z|)~ |f''(z)/f'(z)| \to 0$ as $|z|\to 1-0$, is that the extension (\ref{beta}) provides a quasiconformal extension of $f$, with complex dilatation $\mu$, such that  $\mu(z) \to 0$ as $|z|\to 1+0$. It is known  that the converse holds as well, meaning that if such a quasiconformal extension exists (not necessarily the one given in (\ref{beta})), then 	$\log{f'}\in \mathcal{B}_0$.

Given a space $X$  of analytic functions in  $\D$, we say that a Jordan curve $\Gamma$ is a $X$-curve if the logarithmic derivative of the conformal mapping $f$ sending $\D$ onto the interior of $\Gamma$ is in $X$. In the second part of this note, we will study the geometric properties of the curves $\Gamma$ associated to the spaces of functions we have been considering.

Classical results \cite{p2} state that the $ \mathcal{B}_0$-curves  are the so called  asymptotically conformal curves, which are characterized as follows:
$$
\sup_{z\in \Gamma(z_1,z_2)}  \frac{|z_1-z|+|z-z_2|}{|z_1-z_2|}\rightarrow 1 \quad \operatorname{as}~ |z_1-z_2|\rightarrow 0,
$$
where ${\Gamma}(z_1,z_2)$ denotes the subarc of  $\Gamma$ joining $z_1$ and $z_2$ with smaller diameter. In the $VMOA$ setting, it holds that $\Gamma$ is a  $VMOA$-curve if and only if $\Gamma$ is an asymptotically smooth curve, that is  for all $z_1,z_2 \in \Gamma$
$$
\frac{\ell_{\Gamma}(z_1,z_2)}{|z_1-z_2|} \rightarrow 1 \quad \operatorname{as} ~ {|z_1-z_2|} \rightarrow 0,
$$ 
where $\ell_\Gamma(z_1,z_2)$ is the length along the curve  of $\Gamma(z_1,z_2)$.
We refer the reader to (\cite{bp}, \cite{d}, \cite{p1}, \cite{p2}), where the history and different proofs of these results can be found.\\

The first attempt to geometrically  characterize  Dirichlet-curves, known as Weil-Petterson curves (WP), appeared in \cite{g}, where the authors state that if $\Gamma$ is a WP curve, then 
\begin{equation} \label{ga}
\int_{\Gamma} \int_{\Gamma} \frac{\ell_{\Gamma}(z_1,z_2)-|z_1-z_2|}{|z_1-z_2|^3}~ |dz_1||dz_2|< \infty.
\end{equation}
 Although there is an  error in the proof given in \cite{g},  Ch. Bishop  proves in  \cite{b} that  the  statement  is correct. Moreover, among many other results, he shows, using techniques that involve the beta numbers,  that the expression  in (\ref{ga})   characterizes WP curves, .\\
The beta numbers were introduced by P. Jones in order to study  the rectifiability properties of sets $E\subset \mathbb{R}^2$. Given a dyadic square $Q\subset \mathbb{R}^2$, the beta number $\beta_E(Q)$ associated to the set $E$ is defined as
$$
\beta(Q) = \beta_E(Q) =\frac{1}{\diam(Q)}\inf_L \sup\{\textrm{dist}(z, L) : z \in 3Q \cap E\},
$$
where the infimum is over all lines $L$ that hit $ 3Q.$\\
Ch. Bishop proves in \cite{b2}  a new version of the traveling salesman theorem, which states that for any Jordan arc with endpoints $z_1$ and $z_2$,
\begin{equation}\label{TSP}
\ell_{\Gamma}(z_1,z_2)-|z_1-z_2|\simeq \sum_{Q} \beta^2_\Gamma(Q) \diam(Q),
\end{equation}
where the sum is taken over all dyadic squares $Q\subset \mathbb{R}^2$. As a consequence, it is shown that the condition $\sum_{Q } \beta^2_\Gamma(Q) <\infty 
$,  is equivalent to (\ref{ga}), and therefore  provides a  characterization of  WP curves in terms of the $\beta$-numbers.\\
\indent
For $Q_{p,0}$ spaces $0<p<1$, the only result involving the geometry of the curves known to the author, is due to J. Pau and J.A. Pel\'aez  (\cite{PP}, Th. 4), where they show that  $\Gamma$ is a $Q_{p,0}$-curve if
$$
\int_0^1 \frac{\epsilon^2(t)}{t^{2-p}} dt< \infty,
	$$
	where $\epsilon(t)=\sup_{|z_1-z_2|\leq t} \sup_{z\in \Gamma(z_1,z_2)}  \left(\frac{|z_1-z|+|z-z_2|}{|z_1-z_2|}-1\right)^{1/2}$.
	
We present in the next theorem a result on the other direction. For that, we need to consider a different  version of the beta numbers. Let $\Gamma$ be a Jordan curve, and $\gamma\subset  \Gamma$ be an arc in $\Gamma$ with endpoints $w_1, w_2$. We define
\begin{equation}\label{def}\beta (\gamma)= \frac{1}{|w_1-w_2|} \max \{\textrm{dist}(w, L) : w \in \gamma\},\end{equation}
where $L$ is the line passing through $w_1$  and $w_2$.\\ 
\begin{theo}
	 \begin{itemize}
	 	\item[(i)]
	If a Jordan curve $\Gamma$ is a	
 	$ Q_{p,0}$-curve,~$0<p\leq 1$, then for any dyadic arc $J\subset \T$, we have
 	$$
 	\frac{1}{|J|^p}\sum_{ I\subset J} \beta^2(\gamma_I) ~|I|^p \to 0 ~\textrm{as} ~|J|\to 0,
 	$$
 	where the sum is over all the dyadic intervals $I$ in $\T$ such that $I\subset J$, and $\gamma_I=f(I)$.
 	\item[(ii)] If a Jordan curve $\Gamma$ is a	
 	$ \mathcal{D}_{\log_p}$-curve,~$p\geq 0$, then 
 	$$
 	\sum_{ I} \beta^2(\gamma_I) \left(\log \frac{1}{|I|}\right)^p<\infty,
 	$$
 	where the sum is over all the dyadic intervals $I$ in $\T$.
 	\end{itemize}
 \end{theo}

Next, for any  arc $I \in \T$, with endpoints $z_1,z_2$, we define the quantity  
$$\Delta(I)=\frac{ \ell_\Gamma(w_1,w_2)-|w_1-w_2|}{|w_1-w_2|}$$
where  $w_i=f(z_i), i=1,2$, that is, the image of $z_1, z_2$ under the conformal map $f$ that sends the unit disc onto the interior of $\Gamma$.
Applying  the previous theorem and (\ref{TSP}), we obtain the following corollary:

\begin{coro}\begin{itemize}
		\item[(i)] If  $\Gamma$ is a	
	$ Q_{p,0}$-curve,~$0<p< 1$,   then for any dyadic arc $J\in \T$
	$$
	\frac{1}{|J|^p} \sum_{I} \Delta (I) |I|^p  \to 0~ \textrm{as} ~|J|\to 0,
	$$
	where the sum is over all the dyadic intervals $I\subset J$.\\
		\item[(ii)] 
	Analogously, if  $\Gamma$ is a	
	$ \mathcal{D}_{\log_p}$-curve,~$p\geq 0$, then 
	$$
	\sum_{ I} \Delta(I) \left(\log \frac{1}{|I|}\right)^p<\infty,
	$$
	where the sum is over all the dyadic intervals $I$ in $\T$.
\end{itemize}
	\end{coro}
	Observe that in (i), we do not obtain the result for $p=1$. \\
	These expressions are  the corresponding discrete  versions  of (\ref{ga}) for the spaces we are considering. Note that in this case, the sum is over dyadic arcs on $\T$, or equivalently, over a dyadic partition  on the curve $\Gamma$, but   with respect to   harmonic measure  instead of arclength. We will comment more on this at the end of section 2.

 Throughout this paper, the letter $c$ denotes a constant that may change at different ocurrences. The notation $A \lesssim B~ (A\gtrsim B)$ means that there is a constant $c$ such that $ A \leq cA ~ (A\geq c B)$. Also, as usual, we denote by $\T$ the boundary of the unit disc, and by $B(z, R)$ the ball of radius $R$ centered at the point $z\in \C$.
 
 Finally, the author would like to thank the referee for the very useful comments and remarks that have improved the final results and presentation of this paper, as well as Ch. Bishop and I. Uriarte-Tuero for many helpful conversations.

\section{Quasiconformal characterizations}
\begin{proof}[Proof of the Remark]
	We begin this section by showing that if $\log{f'}\in \mathcal{D}_{\log_p},$ for $p>2$, then $f'$ is Dini smooth in $\overline\D$.\\
	For any point $z\in\D$ with $1-|z|<1/10$, consider the ball  $B_z=B\left(z, \frac{1}{2} (1-|z|) \right)$. Note that if $w\in B_z$, then $\log \left(\frac{1}{1-|w|}\right)>\log \left(\frac{2}{3(1-|z|)}\right)>\frac{1}{2} \log \left(\frac{1}{(1-|z|)}\right)$.\\
	Since $\left|\frac{f''}{f'}\right|$ is subharmonic, for any $z\in \D$ such that $1-|z|<1/10$, it holds that 
	\begin{align}\nonumber
	\left\vert\frac{f''(z)}{f'(z)}\right\vert^2 &\lesssim \frac{1}{(1-|z|)^2} \int_{B_z}\left \vert\frac{f''(w)}{f'(w)}\right\vert^2 ~dA(w)\nonumber\\&\lesssim  \frac{\left(\log \frac{1}{1-|z|}\right)^{-p}}{(1-|z|)^2} \int_{B_z}\left\vert\frac{f''(w)}{f'(w)}\right \vert^2 \left(\log \frac{1}{1-|w|}\right)^p dA(w).
	\end{align}
	The last integral is bounded by the integral in $\D$. Therefore, if $1-|z|<1/10$,
	\begin{equation}\label{dini}
	\left\vert\frac{f''(z)}{f'(z)}\right\vert\lesssim \frac{1}{(1-|z|) \left(\log \frac{1}{1-|z|}\right)^{p/2}}.
	\end{equation}
	This implies that when $p>2$, the integral 
	$$
	\int_{9/10}^1  \left\vert\frac{f''(rz)}{f'(rz)}\right \vert ~dr \lesssim \int_{9/10}^1 \frac{1}{(1-r) \left(\log (1/(1-r)\right)^{p/2}}~dr <\infty
	$$
	converges uniformly in $z\in \T$, and therefore $\log f'$ is continuous in $\overline \D$. In particular, both  $f'$ and $1/f'$ are  continuous in $\overline \D$. Thus, by (\ref{dini}) we have that $|f''(z)|\lesssim \frac{1}{(1-|z|)\left(\log \frac{1}{1-|z|}\right)^{p/2}} $, which by well-known  classical results implies that $f'$ is Dini-continuous in $\T$.

	\end{proof}
Next, let us recall the statement in Theorem \ref{theo 1}: We have to show the equivalence of the conditions 	$\log{f'}\in \mathcal{D}_{\log_p}$ and 
	\begin{equation}\label{D2}
\int_{\C\backslash \D }~ \frac{|\mu(z)|^2}{(|z|^2-1)^2 }\left(\log \frac{|z|}{|z|-1}\right)^p~ dA(z)<\infty, ~p\geq 0.
\end{equation}

\begin{proof}[Proof of Theorem 1]

	Firstly, let us assume that $\log{f'}\in \mathcal{D}_{\log_p},~p\geq 0$, that is,$$ \int_{\D}\left\vert\frac{f''(z)}{f'(z)} \right\vert^2 \left(\log \frac{1}{1-|z|}\right)^p dA(z)<\infty.$$ Since  $ \mathcal{D}_{\log_p} \subset \mathcal{B}_0$, it follows that $ \eta(z)=(1-|z|)~ |f''(z)/|f'(z)| \to 0$ as $|z|\to 1-0$. Therefore,   by the results in  \cite{bp} mentioned in the Introduction, the expression in (\ref{beta}) defines a quasiconformal extension of $f$  to a neighborhood $U$ of $\T$, and from ($\ref{2}$) we deduce that  (\ref{D2}) holds on $U$.\\
	In order to obtain a global quasiconformal mapping, we apply the following theorem (see \cite{lv}, Th.8.1):  If $f:G\to G$' is a $K$-quasiconformal map and $E\subset G$  is a compact set, then there exists a  $K'$-quasiconformal  map in the whole plane that coincides with $f$ in  $E$ , and with $K'$ depending only on $ K, G$ and $E$. \\
	Setting $E=\{z:|z|\leq R_0 \}$, for some $R_0>1$, the above result provides a global quasiconformal extension of the conformal mapping $f$ whose complex dilatation $\mu$ satisfies (\ref{D2}).
	
	To prove the converse, let us assume that (\ref{D2}) holds. Since we are studying the boundary behavior of the map, we can assume without lost of generality that the complex dilation $\mu$ vanishes outside some neighborhood of $\T$, and that  $||\mu||_\infty\neq 0$.\\
	The following  estimate  established by E. Dyn'kin in (\cite{d}, Th. 1) plays a fundamental role in our proof: If $||\mu||_\infty=k<1$
	
	\begin{equation}\label{dyn}
	\bigg\vert\frac{f''(z)}{f'(z)}\bigg\vert \lesssim~ (1-|z|)^{-k} \int_{1-|z|}^{\infty}
	\frac{\omega(z,t)}{t^{2-k}}~ dt, \quad |z| < 1,
	\end{equation}
	
	where
	$$
	\omega^2(z, t) = \frac{1}{\pi t^2}  \int_{|w-z|\leq t} |\mu(w)|^2 ~dA(w).
	$$	
	
	Let $0<a<1$ to be determined later on. By applying the Cauchy-Schwarz inequality, we obtain	
	$$
	\int_{1-|z|}^{\infty}
	\frac{\omega(z,t)}{t^{2-k}}~ dt \leq ~ \left(\int_{1-|z|}^{\infty}
	\frac{\omega^2(z,t)}{t^{2a(2-k)}}~ dt\right)^{1/2} ~	\left(\int_{1-|z|}^{\infty}
	\frac{1}{t^{2(1-a)(2-k)}}~ dt\right)^{1/2}.
	$$
	The right-hand side integral is finite, provided that $2(1-a)(2-k)>1$, i.e.
	$$
	a<\frac{3-2k}{2(2-k)},
	$$	
	in which case
	$$
	\int_{1-|z|}^{\infty}
	\frac{1}{t^{2(1-a)(2-k)}}~ dt\simeq\frac{1}{(1-|z|)^{2(1-a)(2-k)-1}}.
	$$
	\noindent
	These estimates, together with  (\ref{dyn}) and Fubini's theorem, imply that
	\begin{align}
	&\int_{\D} \left\vert\frac{f''(z)}{f'(z)}\right\vert^2 \left( \log \frac{1}{1-|z|}    \right)^p~dA(z)\nonumber\\&\lesssim \int_{\D}\frac{\left( \log \frac{1}{1-|z|}  \right   )^p}{(1-|z|)^{2k+2(1-a)(2-k)-1}}~\left(	 \int_{1-|z|}^{\infty}
	\frac{1}{t^{2a(2-k)+2}} \left(\int_{|w-z|\leq t} |\mu(w)|^2 ~dA(w)\right) dt\right) dA(z) \nonumber\\
	&\lesssim\int_{\C\backslash \D} |\mu(w)|^2\left(\int_{|w|-1}^\infty\frac{1}{t^{2a(2-k)+2}}
	\left(\int_{B(w,t) \cap \D}\frac{\left( \log \frac{1}{1-|z|}\right     )^p}{(1-|z|)^{2k+2(1-a)(2-k)-1}} dA(z)\right)dt\right)dA(w). \label{j}
	\end{align}
	\noindent
	Since $\mu$ has compact support, only the behaviour in a small neighborhood  of $\T$ matters. For $w\in \C\backslash \D $ and $t>|w|-1$ small enough, setting $\tilde{w}=w/|w|$, it holds that $B(w,t)\cap \D \subset B(\tilde{w}, t)\cap \D$. Therefore,  for any $0<s<1$, 
	\begin{equation}\label{h}
	\int_{B(w,t) \cap \D}\frac{1}{(1-|z|)^{s}} \left( \log \frac{1}{1-|z|}     \right)^p~dA(z) \lesssim ~ t^{2-s}( \log (1/t) )^p,
	\end{equation}
	this is because if   $Q=\{(x,y); |x|<l, 0<y<l\},~p\geq 0,~0<s<1$, a simple calculation yields,
$$\int_Q \frac{1}{y^{s}} (\log (1/y))^p dx dy\simeq l^{2-s} (\log (1/l))^p.$$

The condition   $s= 2k+2(1-a)(2-k)-1<1$ holds if and only if $a>1/(2-k)$. Since $k<1/2$, we can choose $0<a<1$ so that $\frac{1}{2-k}<a<\frac{3-2k}{2(2-k)}$.
	Thus, we conclude by (\ref{j}) and  (\ref{h}) that 
	$$
	\int_{\D}	\bigg\vert\frac{f''(z)}{f'(z)}\bigg\vert^2 \left( \log \frac{1}{1-|z|}    \right)^p ~dA(z)\lesssim
	\int_{\C\backslash \D} \frac{|\mu(w)|^2}{(|w|-1)^2} \left( \log \frac{1}{|w|-1}    \right)^p dA(w),
	$$
	which yields the conclusion, since $\mu$ vanishes outside a small neighborhood of $\T$.
	\end{proof}

In order to  prove Theorem 2, we need to show the equivalence of the conditions $\log f'\in Q_{p,0},~0<p\leq 1$, and the dilatation coefficient $\mu$ satisfying:
\begin{equation}\label{qp}
\lim_{I\to 0} \frac{1}{|I|^p}	\int_{\tilde{Q}_I} \frac{|\mu(z)|^2}{(|z|-1)^{2-p} }~ dA(z)=0.
\end{equation}

\begin{proof}[Proof of Theorem 2]
	
Since $Q_{p,0}\subset \mathcal{B}_0$, we proceed as in the previous case to show that if $\log f'\in  Q_{p,0},~0<p\leq 1$, then the expression in (\ref{beta}) provides  a global quasiconformal extension of the conformal mapping $f$ with $\mu$ as in (\ref{qp}).

To prove the converse, we will again apply  the estimate in (\ref{dyn}).

Writing $t^{2-k}=t^{a(2-k)}~t^{(1-a)(2-k)}$ for some $0<a<1$ to be fixed later on, and applying the Cauchy-Schwarz inequality, we obtain	
\begin{equation}\label{ss}
\int_{1-|z|}^{\infty}
\frac{\omega(z,t)}{t^{2-k}}~ dt \leq ~ \left(\int_{1-|z|}^{\infty}
\frac{\omega^2(z,t)}{t^{2a(2-k)}}~ dt\right)^{1/2} ~	\left(\int_{1-|z|}^{\infty}
\frac{1}{t^{2(1-a)(2-k)}}~ dt\right)^{1/2}.
\end{equation}	
By choosing $a$ so that $2(1-a)(2-k)>1$, that is, $	a<\frac{3-2k}{2(2-k)}$, we obtain
\begin{equation}\label{sss}
\int_{1-|z|}^{\infty}\frac{1}{t^{2(1-a)(2-k)}}~ dt \simeq  \frac{1}{(1-|z|)^{2(1-a)(2-k)-1}}.
\end{equation}
Let $\tilde{Q}= \tilde{Q}_I \subset\C\backslash \D $ be the Carleson box associated to an arc $I\subset \T$. Denote by $Q$ the Carleson box in $\D$ which is symmetric to $\tilde{Q}$, that is  $Q=  \{r e^{it}: e^{it}\in I, 1-|I| < r \leq 1\}$.\\
By (\ref{dyn}), (\ref{ss}), (\ref{sss}), and Fubini's theorem,
\begin{align}
&\int_{Q} \bigg\vert\frac{f''(z)}{f'(z)}\bigg\vert^2 (1-|z|)^p ~dA(z)\nonumber\\&\lesssim~ \int_{Q}\frac{1}{(1-|z|)^{2k+2(1-a)(2-k)-1-p}}~	 \left(\int_{1-|z|}^{\infty}
\frac{1}{t^{2+2a(2-k)}}\left( \int_{|w-z|\leq t} |\mu(w)|^2dA(w)\right) dt\right) ~dA(z) \nonumber\\
&\lesssim ~\int_{\C\backslash \D} |\mu(w)|^2\left(\int_{ \{(z,t):z\in Q, |w-z|\leq t\} }\frac{1}{(1-|z|)^{2k+2(1-a)(2-k)-1-p}}~             \frac{1}{t^{2+2a(2-k)}}dt~ dA(z)\right)dA(w)\label{n}\\ 
&=\int_{\{w\in  \tilde Q_N\}}+ \int_{\{w\notin \tilde Q_N \}}=I_1+I_2,\label{nn}
\end{align}
where     $\tilde{Q_n}= 2^n \tilde{Q}; n\geq 1 $. The letter $ N\geq 2$ denotes an integer that will be fixed later on.\\
\indent
If $w\notin \tilde{Q}_N$, and $z\in Q$, we have $|w-z|\simeq |w-\xi_I|$, where $\xi_I$ is the middle point of $I$. Set $s= 2k+2(1-a)(2-k)-1-p$. Provided that $s<1$, by integrating on $t$, the inner integral in (\ref{n}) can be expressed as
$$
\int_ {Q}\frac{1}{(1-|z|)^s}~ \frac{1}{|w-z|^{1+2a(2-k)}}~ dA(z)\simeq  \frac{1}{|w-\xi_I|^{1+2a(2-k)}} ~|I|^{2-s}.
$$
Therefore, 
\begin{align}
I_2&\lesssim |I|^{2-s}~\int_{w\notin \tilde{Q}_N} \frac{|\mu(w)|^2}{|w-\xi_I|^{1+2a(2-k)}}~dA(w)\nonumber\\&= |I|^{2-s}\sum_{n\geq N} \int_{\tilde{Q}_{n+1}\backslash \tilde{Q}_{n}} \frac{|\mu(w)|^2}{|w-\xi_I|^{1+2a(2-k)}}~dA(w)\nonumber\\
&\lesssim |I|^{2-s}\sum_{n\geq N}\frac{1}{(2^n|I|)^{1+2a(2-k)}} \int_{\tilde{Q}_{n+1}}|\mu(w)|^2 dA(w) \nonumber\\&\lesssim \frac{1}{|I|^{2-p}}
\sum_{n\geq N}\frac{(2^n |I|)^{2-p}}{(2^n)^{1+2a(2-k)}} \int_{\tilde{Q}_{n+1}}\frac{|\mu(w)|^2}{(|w|-1)^{2-p}} dA(w)\nonumber\\
&\simeq \sum_{n\geq N} \frac{1}{(2^n)^{2a(2-k)-1+p}}\int_{\tilde{Q}_{n+1}}\frac{|\mu(w)|^2}{(|w|-1)^{2-p}} dA(w).\nonumber
\end{align}
Since the last integral is bounded by $ c~(\diam (\tilde{Q}_{n+1}))^p\simeq (2^n |I|)^p$, by imposing the condition  $2a(2-k)-1>0$, we get 
\begin{equation}\label{4}
I_2\lesssim  |I|^p \sum_{n\geq N}\frac{1}{(2^n)^{2a(2-k)-1}}\leq \epsilon |I|^p,
\end{equation}
for any $\epsilon>0$ as close to $0$ as we wish, by choosing some fixed $N$ big enough.

On the other hand, if $w\in \tilde{Q}_N$, the integral $I_1$ in (\ref{nn}) can be written as 
\begin{equation}\label{nnn}
I_1\simeq\int_{\tilde{Q}_N} |\mu(w)|^2\left(\int_{|w|-1}^\infty  \frac{1}{t^{2+2a(2-k)}} \left(\int_{B(w,t)\cap Q}\frac{1}{(1-|z|)^{s}}   dA(z)\right)dt \right)dA(w).
\end{equation}
Since $s<1$,
\begin{equation}\label{l}
\int_{B(w,t) \cap Q}\frac{1}{(1-|z|)^{s}} dA(z) \lesssim ~ t^{2-s}.
\end{equation}
Thus,  by (\ref{nnn}) and (\ref{l}), we conclude that
\begin{equation}\label{I1}
I_1\lesssim \int_{\tilde{Q}_N}|\mu(w)|^2\left(	
\int_{|w|-1}^\infty \frac{1}{t^{3-p}} dt\right)dw \lesssim \int_{\tilde{Q}_N} \frac{|\mu(w)|^2}{(|w|-1)^{2-p}} dA(w)=o(|I|^p),
\end{equation}
since $\mu$ satisfies (\ref{qp}).\\
\indent
The last estimate (\ref{I1}), together with (\ref{4})  and (\ref{nn}), complete the proof that $\log f'\in Q_{p,0}$, as long as there exists $0<a<1$ satisfying $$s= 2k+2(1-a)(2-k)-1-p<1, \quad	a<\frac{3-2k}{2(2-k)}~\textrm{ and}\quad 2a(2-k)-1>0.$$
An immediate computation shows that these conditions hold for some $0<a<1$, provided that $2-p<3-2k$, or equivalently $k<(1+p)/2$.

\end{proof}

 \section{Geometric properties}
 Before proceeding to the proof of Theorem 3, let us introduce some notation. For each dyadic interval $I\subset \T$, we denote by $Q_I\subset\D$ its corresponding Carleson box and by $W_I=T(Q_I)$ the top of $Q_I$, that is $W_I=\{r e^{it}\in Q_I; e^{it}\in I, |I|/2 < 1-r <|I|\}$. Also denote by $\mathcal{W}$ the collection of all $W_I$ associated to the dyadic intervals $I\subset T$.
 For any $W\in \mathcal{W}$, we define $\eta(W)=\sup_{z\in W} |\frac{f''(z)}{f'(z)}| (1-|z|)$.\\
 Let $z_I \in W_I = T(Q_I) \in \mathcal{W}$, and suppose $ z \in 2Q_I$. Let $W_I = W_0, . . . ,W_N$
 be the list of the tops of Carleson  boxes hit by the hyperbolic geodesic from $z_I$ to $ z$. We will call this chain of $W's$ associated to $z, ~C(z)$.
 The following result due to Ch. Bishop, that we state here as a proposition, will be our main tool. The proof can be found in (\cite{b}, sect 4) under the normalization $f'(z_I)=1$. Nevertheless, for the sake of completeness, we will include the details in the Appendix.
 
 \begin{propo}Let $\Gamma$ be an asymptotically conformal curve, and f the conformal map from $\D$ onto the interior of $\Gamma$. Let $z_I \in W_I = T(Q_I) \in \mathcal{W}$, and suppose $ z \in 2Q_I$. Then
 	\begin{equation}\label{chain}
 	f(z)-f(z_I)=f'(z_I)(z-z_I)+O\left(\diam (f(W_I )) \sum_{W\in C(z)}
 	\eta(W) \left(\frac{\diam (W)} {\diam (W_I)}\right)^\alpha\right).
 	\end{equation}
 	for some $0<\alpha<1$. Moreover, $\alpha$ can be as close to 1 as we need by choosing  $\diam (W_I)$ small enough.
 \end{propo}

\begin{proof}[Proof of Theorem 3]
\indent	
	\begin{itemize}
\item[(i)]Since $\Gamma$ is a $Q_{p,0}$- curve, $0<p\leq1$, for a dyadic arc $I\subset\T$, it holds that
\begin{equation}\label{qp0}
		 \lim_{|I|\to 0}\frac{1}{|I|^p}\int_{Q_I}\bigg\vert\frac{f''(z)}{f'(z)}\bigg\vert^2 (1-|z|)^p dA(z)= 0.
		 \end{equation}
On the other hand, note that for any $z\in W\in\mathcal{W}$,  $1-|z|\simeq \diam (W)$. Moreover, due to the subharmonicity of the function $|f''/f'|$, we obtain 
$$
\eta^2(W) \lesssim  \frac{1}{(\diam( W))^2}  \int_{cW}\bigg\vert\frac{f''(z)}{f'(z)}\bigg\vert^2 (1-|z|)^2dA(z)\simeq \int_{cW}\bigg\vert\frac{f''(z)}{f'(z)}\bigg\vert^2 dA(z),
$$
for some constant $c>1$.\\ 
	Therefore, for any dyadic box $Q$, 
\begin {align}\nonumber
  \sum_{W\in  Q} \eta(W)^2 ~(\diam( W))^p &\lesssim  \sum_{W\in  Q}
    \int_{cW}\left\vert\frac{f''(z)}{f'(z)}\right\vert^2 dA(z) (\diam (W))^p\\\nonumber&\simeq 
 \int_{cQ} \bigg\vert\frac{f''(z)}{f'(z)}\bigg\vert^2 (1-|z|)^p  ~dA(z).
\end{align}
The expression (\ref{qp0}) can then be written as
	\begin{equation}\label{v}
	\frac{1}{|I|^p}	\sum_{W\subset  Q_I} \eta(W)^2 ~(\diam(W))^p \to 0~ \textrm{as} ~|I|\to 0.
	\end{equation}
The equation in (\ref{chain}) shows that $f$ is almost linear on the segment $2I$, meaning that the arc of $\Gamma, ~ \gamma_I= f(I)$, deviates from a straight line at most $$\alpha(W_I)=\sup_{z\in 2I}  \left(\diam( f (W_I)) \sum_{W\in C(z) } \eta(W) \left(\frac{\diam(W)}{\diam (W_I)}\right)^\alpha\right).$$ 
Note that by the circular distortion theorem for quasiconformal mappings, and by  Koebe's theorem, $\diam(f(W_I)) \simeq |I| |f'(z_I) | \simeq \diam (f(I))$.\\
In particular, for each dyadic interval $I\in \T$, and  $W_I = T(Q_I ) \in \mathcal{W}$, we can choose a
boundary point $z \in 2I$  so that, the
corresponding sum along its chain $C(z)$ is comparable to $\alpha (W_I)$. We  call this chain $C(W_I)$. \\
Since $\gamma_I=f(I)$ is a quasiarc, $\diam(f(I)) \simeq |w_1-w_2|$, where $w_1, w_2$ are the endpoints of $\gamma_I$. Thus, by the definition of $\beta (\gamma_I)$ given in (\ref{def}) we get
\begin{equation}\label{alpha}
\beta(\gamma_I) \lesssim \frac {\alpha(W_I) }{\diam (\gamma_I)}\lesssim \sum_{W'\in C(W_I)} \eta(W') \left(\frac{\diam(W')}{\diam (W_I)}\right)^\alpha.
\end{equation}
Choose $1/2<\alpha<1$ and $s>0$ so that  $p<s<2\alpha$. By the Cauchy-Schwarz inequality,
$$
\beta(\gamma_I)\leq \left( \sum_{W'\in C(W_I)} \eta^2 (W') \left(\frac{\diam(W')}{\diam (W_I)}\right)^s\right)^{1/2}  \left(\sum_{W'\in C(W_I)}  \left(\frac{\diam(W')}{\diam (W_I)}\right)^{2\alpha-s}~\right)^{1/2}.
$$
 Note that the last sum is a geometric sum, which is finite since $2\alpha-s>0$.	
Therefore, given a dyadic arc $J\subset\T$ of small enough diameter,
\begin{align}\nonumber
\frac{1}{| J|^p}\sum_{ I: I~\textrm{dyadic}  \subset  J} \beta^2 (\gamma_ I)~| I|^p &\lesssim \frac{1}{| J|^p}\sum_{ I: I~\textrm{dyadic}  \subset  J} \sum_{W'\in C(W_I)} \eta^2(W') \left(\frac{\diam(W')}{\diam (W_I)}\right)^s
        |I|^p\\\label{ca}
&\lesssim \frac{1}{| J|^p}
\sum_{W'\subset Q_J}\eta^2 (W') ~(\diam (W'))^s \sum_{\{I: W'\in C(W_I)\}}\frac{1}{ (\diam (W_I))^{s-p} }     
\end{align}
To estimate the last sum, observe that  $W'=T(Q_{I'})$ for some dyadic interval $I'$. If $I$ is such that  $W'\in C(W_I)$, then $I'\subset c I$ for some constant $c>0$. Since $s>p$, we get that the last sum, which is a geometric  sum, can be bounded by a multiple of $\frac{1}{ (\diam (W'))^{s-p}}$. Thus, by (\ref{ca}) and (\ref{v}), we deduce that
\begin{align}\nonumber
\frac{1}{| J|^p}\sum_{ I: I~\textrm{dyadic}  \subset  J} \beta^2 (\gamma_ I)~| I|^p& \lesssim 
\frac{1}{| J|^p}\sum_{W'\subset Q_J}\eta^2 (W') ~(\diam (W'))^p
 \to 0 ~\textrm{as} ~|J|\to 0,
\end{align}
as we wanted to show.
\item[(ii)]
The proof follows a similar argument as in (i). The condition for  $\mathcal{D}_{\log_p}$-curves, $p\geq 0$, can be written now as:
\begin{equation}\label{last}
	\sum_{W\subset  \D} \eta(W)^2 ~\left(\log \left(1/\diam(W)\right)\right)^p <\infty.
\end{equation}
Moreover, applying the Cauchy-Schwarz inequality in (\ref{alpha}), we get
$$
\beta(\gamma_I)\leq \left( \sum_{W'\in C(W_I)} \eta^2 (W') \left(\frac{\diam(W')}{\diam (W_I)}\right)^\alpha\right)^{1/2}  \left(\sum_{W'\in C(W_I)}  \left(\frac{\diam(W')}{\diam (W_I)}\right)^{\alpha}~\right)^{1/2}_.
$$ 
Since $\alpha>0$, the last sum converges.
Thus, 
\begin{align}\nonumber
\sum_{ I: I~\textrm{dyadic}  \subset  \T} \beta^2 (\gamma_ I)&~(\log (1/|I|))^p\lesssim \sum_{ I: I~\textrm{dyadic}  \subset  \T} 
\sum_{W'\in C(W_I)} \eta^2(W') \left(\frac{\diam(W')}{\diam (W_I)}\right)^\alpha (\log (1/|I|))^p\\\nonumber
&\lesssim 
\sum_{W'\subset \D}\eta^2 (W') ~(\diam (W'))^\alpha \sum_{\{I: W'\in C(W_I)\}}\frac{1}{ (\diam (W_I))^{\alpha} }  (\log (1/|I|))^p. 
\end{align}
Since $0<\alpha<1, ~p\geq 0$, and $\diam (W_I)\simeq |I|$, the last sum is  bounded by a multiple of $\frac{1}{ (\diam (W'))^{\alpha} }  (\log (1/\diam (W')))^p $, which yields
$$
\sum_{ I: I~\textrm{dyadic}  \subset  \T} \beta^2 (\gamma_ I)~(\log (1/|I|))^p\lesssim
\sum_{W'\subset \D}\eta^2 (W') (\log (1/\diam (W')))^p<\infty.$$
because of (\ref{last}).

\textit{Remark 1}:
Since $f$ is quasiconformal in $\C$, there is  $0<\alpha<1$, so that $\left(\frac{|I|}{\diam(\T|)}\right)^{1/\alpha}\leq \frac{\diam{f(I)}}{\diam (f(\T))}\leq \left(\frac{|I|}{\diam(\T|)}\right)^{\alpha}$ (see \cite{AIM}, corollary 3.10.4).\\
Therefore, we obtain  that $\sum_{ I: I~\textrm{dyadic}  \subset  \T} \beta^2 (\gamma_ I)~(\log (1/|I|))^p<\infty$ is equivalent to 
\begin{equation}\label{o}
\sum_{ I: I~\textrm{dyadic}  \subset  \T} \beta^2 (\gamma_ I)~\left(\log \left(\frac {1}{\diam(\gamma_I)}\right)\right)^p<\infty.
\end{equation}

\end{itemize}
\end{proof}

\textit{Remark 2}: We will assume that the reader is familiar with the concept of multi-resolution families, and the relations between the different definitions of beta-numbers (see, Appendix B in  \cite{b2} for an excellent exposition).\\
 It is is obvious that the results in Theorem 3 still hold if we consider any translation of the dyadic intervals, in particular the 1/3-translation. With this "1/3-trick", one generates a multi-resolution family  on $\T$ which is sent by the conformal map to a family of arcs in  $\Gamma$. Let us denote this family  by $\mathcal{\tilde{F}}$. It turns out that $\mathcal{\tilde{F}}$ is a MR family when $\Gamma$ is a chord-arc curve. Applying now Lemma B.1 in  \cite{b2} to the expression (\ref{o}), we obtain that:

If $\log f'\in \mathcal{D}_{\log_p},~p\geq 0$, then
 
$$\sum_{ Q} \beta^2_ \Gamma (Q)~\left(\log \left(\frac {1}{\diam(Q)}\right)\right)^p<\infty,$$
where the sum is taken over all dyadic squares $Q \subset \mathbb{R}^2$.

\begin{proof}[Proof of Corollary 1]
	
 $(i)$ In \cite{b2}, Ch. Bishop proves and discusses  equivalent formulations of the traveling salesman theorem expressed in terms of different multi-resolution families, (see, Appendix B in  \cite{b2}). In our setting, given a dyadic arc $I\subset \T$ with endpoints $z_1, z_2$  and $w_i=f(z_i), i=1,2$, the traveling salesman theorem states that
$$\ell_\Gamma(w_1,w_2)-|w_1-w_2| \simeq \sum_{I' \subset \mathcal{F}_I } \beta^2(\gamma_{I'} )\diam( \gamma_{I'}),$$
where $\gamma_{I'}=f(I')$, and $\mathcal{F}_I$ is the MR family generated by the 1/3 translates of the dyadic intervals contained in $I$.\\
Recall that by definition,
$$
\Delta (I )= \frac{\ell_\Gamma(w_1,w_2)-|w_1-w_2|}{|w_1-w_2|}.
$$
Since  $\diam( f(I))\simeq |w_1-w_2|$,  we get that 
$$\frac{1}{|J|^p}\sum_{I \textrm{dyadic}\subset J} \Delta (I) |I|^p~	\simeq \frac{1}{|J|^p}\sum_{I \textrm{dyadic}\subset J} 
\sum_{I'  \subset \mathcal{F}_I} \beta^2 (\gamma_I')\frac{\diam( f(I'))}{\diam( f(I))}  |I|^p.
	$$
	Since the curve $\Gamma$ is  asymptotically conformal, there   is $0<\alpha<1$ as close to 1 as we wish, such that for any arc $I'\subset I \subset \T$,
	$$
	\frac{\diam( f(I'))}{\diam( f(I))}\lesssim \left(\frac{|I'|}{|I|}\right)^\alpha.
	$$
	Therefore, 
	$$\frac{1}{|J|^p}\sum_{I \textrm{dyadic}\subset J} \Delta (I) |I|^p\lesssim \frac{1}{|J|^p} \sum_{I'  \subset \mathcal{F}_J} \beta^2 (\gamma_{I'})|I'|^\alpha \sum_{I:I'  \subset \mathcal{F}_I  }|I|^{p-\alpha}.
	$$
	The last sum is a geometric sum. Choosing  $1>\alpha>p$, this geometric sum  is bounded by a multiple of $|I'|^{p-\alpha}.$ Hence, applying Theorem 3, we get
	$$\frac{1}{|J|^p}\sum_{I \textrm{dyadic}\subset J} \Delta (I) |I|^p\lesssim \frac{1}{|J|^p} \sum_{I'  \subset \mathcal{F}_J} \beta^2 (\gamma_{I'})|I'|^p\to 0~\textrm{as}\ |J|\to 0.$$
	
	$(ii)$ Proceeding as in the previous case, we get
	\begin{align}\nonumber
	\sum_{I \textrm{dyadic}\subset \T} \Delta (I) (\log (1/ |I|))^p& \lesssim \sum_{I \textrm{dyadic}\subset \T} 
	\sum_{I'  \subset \mathcal{F}_I} \beta^2 (\gamma_{I'}) \left(\frac{|I'|}{|I|}\right)^\alpha (\log (1/ |I|))^p\\\nonumber&
		\lesssim \sum_{I'  \subset \mathcal{F}_\T}\beta^2 (\gamma_{I'})  |I'|^\alpha \sum_{I:I'  \subset \mathcal{F}_I}\frac{1}{|I|^\alpha}(\log (1/ |I|))^p\\\nonumber&\lesssim \sum_{I' \textrm{dyadic}\subset \T}\beta^2 (\gamma_{I'}) (\log (1/ |I'|))^p<\infty,\nonumber
		\end{align}
		because of  Theorem 3 (ii).\\

\end{proof}

\section{Appendix}

In this section we include a sketch of the proof of Proposition 1. We will follow the argument in \cite{b}.\\
Suppose $z_I \in W_I = T(Q_I) \in \mathcal{W}$, and suppose $ z \in 2Q_I$. Let $W_I = W_0, . . . ,W_N$
be the list of the tops of Carleson boxes hit by the hyperbolic geodesic from $z_I$ to $ z$. Then
$$
|\log f'(z)- \log f'(z_I)|=\left|\int_{z_I}^z \frac{f''(\xi)}{f'(\xi)}d\xi\right|\lesssim \sum_{n=0}^{N}
\eta(W_n).
$$
Thus,
\begin{equation}\label{logest}
\left|\frac{f'(z)}{f'(z_I)}-1~\right|=\left|e^{\log f'(z)- \log f'(z_I)}-1\right|\leq e^{|\log f'(z)- \log f'(z_I)|}-1\leq e^{c \sum_{n=0}^{N}
	\eta(W_n)}-1.
\end{equation}
Let us assume first that $\sum_{n=0}^{N}
\eta(W_n)\leq 1$, and let $0<\alpha<1$. Since $e^x-1\lesssim x$ for $x\in[0,1]$, we obtain by (\ref{logest}) that
\begin{align}\nonumber
|f(z)-f(z_I)&-(z-z_I) f'(z_I)|=\left|\int_{z_I}^z f'(\xi)-f'(z_I)~ d\xi\right|\leq  \sum_k\int_{z_k}^{z_{k+1}} |f'(\xi)-f'(z_I)|~ |d\xi|\\\nonumber
&\lesssim |f'(z_I)| \sum_k\diam(W_k) \sum_{n=0}^{k}\eta(W_n)
\lesssim |f'(z_I)| \sum_n \eta(W_n) \sum_{k\geq n} \diam(W_k)\\\nonumber
&\lesssim |f'(z_I)|\sum_n \eta(W_n) \diam(W_n) \lesssim |f'(z_I)|\diam(W_0) \sum_n \eta (W_n) \left(\frac {\diam (W_n)}{\diam (W_0)}\right)^\alpha \\\nonumber
&\simeq  \diam (f(W_I )) \sum_{W\in C(z)}
\eta(W) \left(\frac{\diam (W)} {\diam (W_I)}\right)^\alpha.\label{es}\\
\end{align}
 We have used that fact that the $\diam(W_k)$ decrease geometrically, and
hence the sum is dominated by a multiple of its largest term.

Assume next that $\sum_{n=0}^{N}\eta(W_n)> 1$. We can then choose $M$ along the hyperbolic geodesic such that  $\sum_{n=0}^{M} \eta(W_n)\leq 1$ and   $\sum_{n=0}^{M+1} \eta(W_n)\geq 1$. On one hand, by (\ref{es}),
\begin{equation}\label{a1}
|f(z_M)-f(z_I)-(z_M-z_I) f'(z_I)|\lesssim \diam (f(W_I )) \sum_{n=0}^{M}
\eta(W_n) \left(\frac{\diam (W_n)} {\diam (W_0)}\right)^\alpha.
\end{equation}
On the other hand, since $\log f'\in \mathcal{B}_0$, for all
sufficiently small Carleson boxes, the conformal map
$f$ restricted to such a Carleson box $2Q_I$ has a $K$-quasiconformal extension to the
reflection across $2I$, with $K$ close to 1. By Mori's theorem we have that 
$$
|f(z)-f(z_M)|\leq C~|f(z)-f(z_I)|~\bigg|\frac{z-z_M}{z-z_I}\bigg|^\alpha.
$$
for some $C = C(\alpha) < \infty$, where we may take $\alpha < 1$ as close to $1$ as we wish.\\
Therefore,
\begin{align}\nonumber
|f(z)-f(z_M)|&\lesssim \diam (f(W_0))\left(\frac{\diam(W_M)}{\diam (W_0)}\right)^\alpha\\\nonumber
&\leq \diam (f(W_0))\left( \sum_{n=0}^{M} \eta(W_n)\right) \left(\frac{\diam(W_M)}{\diam (W_0)}\right)^\alpha\\& \lesssim \diam (f(W_0)) \sum_{n=0}^{M} \eta(W_n) \left(\frac{\diam(W_n)}{\diam (W_0)}\right)^\alpha,\label{a2}
\end{align}
since $\diam (W_n)\geq \diam(W_M)$ for all $n\leq M$. Similarly,
\begin{align}\nonumber
|z-z_M|&=|z-z_I|~\bigg|\frac{z-z_M}{z-z_I}\bigg|\leq |z-z_I|~\bigg|\frac{z-z_M}{z-z_I}\bigg|^\alpha \\
&\leq \diam (W_0) \sum_{n=0}^{M} \eta(W_n) \left(\frac{\diam(W_n)}{\diam (W_0)}\right)^\alpha.\nonumber
\end{align}
By Koebe's theorem,
\begin{equation}\label{a3}
|(z-z_M) f'(z_I)|\lesssim \diam \left(f (W_0)\right) \sum_{n=0}^{M} \eta(W_n) \left(\frac{\diam(W_n)}{\diam (W_0)}\right)^\alpha.
\end{equation}
The result follows from (\ref{a1}), (\ref{a2}) and (\ref{a3}), since
\begin{align}\nonumber
|f(z)-f(z_I)&-(z-z_I) f'(z_I)|\\\nonumber
&\leq 
|f(z)-f(z_M)-(z-z_M) f'(z_I)|+ |f(z_M)-f(z_I)-(z_M-z_I) f'(z_I)|\\\nonumber
&	\leq 	|f(z)-f(z_M)|+|(z-z_M) f'(z_I)|+ |f(z_M)-f(z_I)-(z_M-z_I) f'(z_I)|\\
&\lesssim  \diam (f (W_0)) \sum_{n=0}^{M} \eta(W_n) \left(\frac{\diam(W_n)}{\diam ((W_0))}\right)^\alpha\nonumber\\&\lesssim  \diam (f(W_I )) \sum_{W\in C(z)}
\eta(W) \left(\frac{\diam (W)} {\diam (W_I)}\right)^\alpha.\nonumber
\end{align}

\textit{Mar\'ia J. Gonz\'alez:} Departamento de Matem\'aticas, Universidad de C\'adiz, 11510 Puerto Real (C\'adiz), Spain. E-mail address: majose.gonzalez@uca.es

\end{document}